\documentclass[a4paper]{amsart}
\usepackage{amssymb}
\usepackage{amsaddr}
\usepackage{xypic}
\usepackage[only,mapsfrom]{stmaryrd}
\usepackage{graphicx}

\newtheorem{theorem}{Theorem}[section]
\newtheorem{lemma}[theorem]{Lemma}

\newtheorem{cor}[theorem]{Corollary}

\theoremstyle{definition}
\newtheorem{definition}[theorem]{Definition}

\theoremstyle{remark}
\newtheorem{remark}[theorem]{Remark}

\numberwithin{equation}{section}

\newcommand{\C}{\mathbb{C}}

\newcommand{\R}{\mathbb{R}}

\DeclareMathOperator{\image}{im}

\DeclareMathOperator{\id}{Id}
\DeclareMathOperator{\aut}{aut}
\DeclareMathOperator{\Diff}{Diff}
\DeclareMathOperator{\Homeo}{Homeo}

\title{Non-abelian symmetries of quasitoric manifolds}

\author{Michael Wiemeler}
\address{Institut f\"ur Algebra und Geometrie\\ Karlsruher Institut f\"ur Technologie\\ Kaiserstrasse 89-93\\ D-76133 Karlsruhe\\ Germany}
\email{michael.wiemeler@kit.edu}
\thanks{Parts of the research were supported by SNF Grant No. PBFRP2-133466 and a grant from  the MPG}



\begin{document}
\begin{abstract}
  A quasitoric manifold \(M\) is a \(2n\)-dimensional manifold which admits an action of an \(n\)-dimensional torus which has some nice properties.
  We determine the isomorphism type of a maximal compact connected Lie subgroup \(G\) of \(\Homeo(M)\) which contains the torus.
  Moreover, we show that this group is unique up to conjugation.
\end{abstract}

\maketitle

%-----------------------------------------------------------------------
% End of amsart.template
%-----------------------------------------------------------------------

\section{Introduction}
\label{sec:intro}

A quasitoric manifold is a smooth connected orientable \(2n\)-dimensional manifold \(M\) with a smooth action of an \(n\)-dimensional torus \(T\) such that:
\begin{itemize}
\item The \(T\)-action on \(M\) is locally standard, i.e. the \(T\)-action is locally modelled on the standard \(T\)-action on \(\C^n\).
\item If the first property is satisfied, then \(M/T\) is naturally an \(n\)-dimensional manifold with corners.
We require that the orbit space of the \(T\)-action on \(M\) is face-preserving homeomorphic to an \(n\)-dimensional simple polytope.
\end{itemize}

Quasitoric manifolds were introduced by Davis and Januszkiewicz \cite{0733.52006} in 1991.
A symplectic \(2n\)-dimensional manifold with an hamiltonian action of an \(n\)-dimensional torus is an example of a quasitoric manifold.
We call such a manifold a symplectic toric manifold.

McDuff and Tolman \cite{1202.52010} and Masuda \cite{1210.53081} independently constructed a maximal compact connected Lie subgroup of the symplectomorphsim group of a symplectic toric manifold which contains the torus.
Masuda also asked if this maximal Lie subgroup is unique up to conjugation.

In this paper we construct a maximal compact connected Lie subgroup of the homeomorphism group of a quasitoric manifold \(M\) which contains the torus \(T\).
Moreover, we give a partial answer to Masuda's question.
 To be more precise, we have the following theorem.

\begin{theorem}
\label{sec:introduction}
  Let \(M\) be a quasitoric manifold. Then there is a compact connected Lie subgroup \(G\) of \(\Homeo(M)\) which contains the torus \(T\) such that:
  \begin{enumerate}
  \item \(G\) acts smoothly on \(M\) for some smooth structure on \(M\).
  \item If \(G'\subset \Homeo(M)\) is another compact connected Lie subgroup which contains the torus \(T\) and acts smoothly on \(M\) for some smooth structure on \(M\), then \(G'\) is conjugated in \(\Homeo(M)\) to a subgroup of \(G\).
  \item If the \(G\)- and \(G'\)-actions are smooth with respect to the same smooth structure on \(M\), then \(G'\) is conjugated in \(\Diff(M)\) to a subgroup of \(G\).
  \item If \(M\) is a symplectic toric manifold, then \(G\) is conjugated in \(\Homeo(M)\) to a subgroup of the symplectomorphism group of \(M\).
  \end{enumerate}
\end{theorem}

A smooth structure on \(M\) for which the \(G\)-action from the above theorem is smooth can be described as follows.
By Theorem 5.6 of \cite{wiemeler_pre4}, the \(T\)-equivariant smooth structures on \(M\) correspond one-to-one to smooth structures on the orbit space \(M/T\).
The \(G\)-action is smooth for the \(T\)-equivariant smooth structure on \(M\) for which \(M/T\) is diffeomorphic to a simple polytope.

In this paper all actions of compact Lie groups on manifolds \(M\) are smooth with respect to some smooth structure on \(M\).

This article is organised as follows.
In  Section~\ref{sec:pairs} we review some basic facts about quasitoric manifolds and introduce an automorphism group for the characteristic pair corresponding to a quasitoric manifold.
In Section~\ref{sec:construct} we construct the group \(G\) from Theorem~\ref{sec:introduction}.
In Section~\ref{sec:class} we review the classification of quasitoric manifolds with \(G\)-action up to \(G\)-equivariant diffeomorphism given in \cite{torus}.
Moreover, we give a classification of  these manifolds up to \(G\)-equivariant homeomorphism.
In Section~\ref{sec:unique} we apply the results of the previous section and show that the group \(G\) has the properties described in Theorem~\ref{sec:introduction}.

\section{Characteristic pairs and their automorphism groups}
\label{sec:pairs}

Let \(M\) be a \(2n\)-dimensional quasitoric manifold, \(P\) its orbit polytope and \(\pi:M\rightarrow P\) the orbit map.
Denote by \(\mathfrak{F}(M)\) the set of facets of \(P\).
We write also \(\mathfrak{F}\) instead of \(\mathfrak{F}(M)\) if it is clear from the context which quasitoric manifold is meant.
Then the preimage \(M_i=\pi^{-1}(F_i)\) of \(F_i\in \mathfrak{F}\) is a codimension two submanifold of \(M\) which is fixed by a one-dimensional subtorus \(\lambda(F_i)=\lambda(M_i)\) of \(T\).
These \(M_i\) are called characteristic submanifolds of \(M\).
Since the facets of \(P\) correspond one-to-one to the characteristic submanifolds of \(M\), we denote the set of characteristic manifolds also by \(\mathfrak{F}\).

Let \(IT\subset LT\) be the integral lattice in the Lie algebra of \(T\).
The characteristic map \(\lambda:\mathfrak{F}\rightarrow \{\text{one-dimensional subtori of } T\}\) lifts to a map \(\bar{\lambda}:\mathfrak{F}\rightarrow IT\cong \mathbb{Z}^n\) such that, for a subset \(\sigma\) of \(\mathfrak{F}\) with \(\bigcap_{F_i\in \sigma}F_i\neq \emptyset\), \(\{\bar{\lambda}(F_i);\;F_i\in\sigma\}\) is part of a basis of \(IT\).
Note that each \(\bar{\lambda}(M_i)\) is unique up to sign.
We call \(\bar{\lambda}\) a characteristic function for \(M\).

Dual to \(P\) there is a simplicial complex \(K\) with vertex set \(\mathfrak{F}\).
A subset \(\sigma\subset \mathfrak{F}\) is a simplex of \(K\) if and only if \(\bigcap_{F_i\in \sigma}F_i\neq \emptyset\).

Note that \(M\) is determined by the combinatorial type of \(P\) (or \(K\)) and \(\bar{\lambda}\) up to equivariant homeomorphism \cite[Proposition 1.8]{0733.52006}.
This construction motivates the following definition.

\begin{definition}
  Let \(K\) be a simplicial complex of dimension \(n-1\) with vertex set \(\mathfrak{F}\).
  Moreover, let \(T\) be an \(n\)-dimensional torus and \(\bar{\lambda}:\mathfrak{F}\rightarrow IT\cong \mathbb{Z}^n\) a map such that for all simplices \(\sigma\) of \(K\) the set \(\{\bar{\lambda}(F_i);\;F_i\in\sigma\}\) is part of a basis of \(IT\).
Then we call \((K,\bar{\lambda})\) a characteristic pair.
\end{definition}

An omniorientation of a quasitoric manifold \(M\) is a choice of orientations for \(M\) and all characteristic submanifolds of \(M\).
An omniorientation of \(M\) induces a complex structure on all normal bundles of the characteristic submanifolds.
These complex structures may be used to make the map \(\bar{\lambda}\) unique by requiring that the \(S^1\)-action induced by \(\bar{\lambda}(M_i)\) on the normal bundle of \(M_i\) is given by complex multiplication.

The cohomology \(H^*(M;\mathbb{Z})\) was computed by Davis and Januszkiewicz \cite[Theorem 4.14]{0733.52006}.
It is torsion-free and generated by the Poincar\'e duals \(PD(M_i)\) of the characteristic manifolds.
If we choose \(\bar{\lambda}\) as above, then these Poincar\'e duals are subject to the following relations
\begin{equation}
\label{eq:1}
  0=\sum_{M_i\in \mathfrak{F}} \langle v,\bar{\lambda}(M_i)\rangle PD(M_i)\quad \text{ for all } v \in IT^*
\end{equation}
and, for \(\sigma\subset \mathfrak{F}\),
\begin{equation*}
  0=\prod_{M_i\in \sigma}PD(M_i) \quad \Leftrightarrow \quad \sigma \text{ is not a simplex of } K.
\end{equation*}

If \(\sigma\) is a simplex of \(K\), then, because \(\{\bar{\lambda}(F_i);\;F_i\in\sigma\}\) is part of a basis of \(IT\), we may choose an isomorphism \(IT\rightarrow \mathbb{Z}^n\), such that the \(\bar{\lambda}(F_i)\)'s are given by the columns of a matrix of the form
\begin{equation*}
  \Lambda=\left(
    \begin{matrix}
      1&&\\
      &\ddots&&\Lambda'\\
      &&1
    \end{matrix}
\right),
\end{equation*}
where the first \(\#\sigma\) columns of \(\Lambda\) correspond to the \(F_i\in \sigma\).
We call such a matrix a characteristic matrix of \(M\).
If \(v_1,\dots,v_n\) is the basis of \(IT^*\) dual to the standard basis of \(IT\cong \mathbb{Z}^n\), then the coefficient  \(\langle v_j, \bar{\lambda}(M_i)\rangle\) in equation (\ref{eq:1}) is the \(i\)-th entry of the \(j\)-th row of \(\Lambda\).
Hence, one can read off all relations in (\ref{eq:1}) from the matrix \(\Lambda\).

We are interested in the symmetries of \(M\).
Since quasitoric manifolds are determined by their characteristic pairs, we should also study automorphisms of characteristic pairs.
Therefore we define:

\begin{definition}
  Let \((K,\bar{\lambda})\) be a characteristic pair. Then we define the automorphism group of \((K,\bar{\lambda})\) to be
  \begin{equation*}
    \aut(K,\bar{\lambda})=\{(f,g)\in \aut(K)\times \aut(T);\; Lg \circ \bar{\lambda}=\bar{\lambda}\circ f\}.
  \end{equation*}
\end{definition}

\begin{lemma}
\label{sec:char-pairs-their-2}
  Let \((K,\bar{\lambda})\) be a characteristic pair. Then the natural map
  \(\aut(K,\bar{\lambda})\rightarrow \aut(K)\) is injective.
\end{lemma}
\begin{proof}
  Let \((f,g)\) be an element of the kernel of this map.
  Then we must have \(f=\id\) and \(Lg(\bar{\lambda}(F_i))=\bar{\lambda}(F_i)\) for all \(F_i\in \mathfrak{F}\).
  Therefore the statement follows because \(LT\) is generated by the \(\bar{\lambda}(F_i)\), \(F_i\in \mathfrak{F}\).
\end{proof}

We sometimes have to choose a special type of omniorientation on a quasitoric manifold.
Therefore we make the following definition.
We call an omniorientation on a quasitoric manifold \(M\) strong if for any two characteristic submanifolds \(M_1,M_2\subset M\) we have
  \begin{equation*}
    \text{PD}(M_1)=\pm\text{PD}(M_2)\;\Rightarrow \; \text{PD}(M_1)=\text{PD}(M_2).
  \end{equation*}

The term strong is motivated by the fact that, if the omniorientation on \(M\) satisfies the above condition, then each cone in \(H^2(M;\mathbb{Z})\otimes \R\) which is spanned by two Poincar\'e duals of characteristic manifolds is strongly convex, i.e. this cone does not contain a straight line.

For \(\alpha\in H^2(M)\), we also define 
\begin{equation*}
\mathfrak{F}^{\alpha}=\{ M_i\in \mathfrak{F};\; \text{PD}(M_i)=\alpha\}.  
\end{equation*}

\begin{lemma}
\label{sec:char-pairs-their}
  Let \(M\) be a quasitoric manifold with a strong omniorientation.
 
Then there is a unique homomorphism \(\phi:\prod_{\alpha\in H^2(M)}S(\mathfrak{F}^\alpha)\hookrightarrow \aut(K,\bar{\lambda})\) such that \(\psi\circ\phi:\prod_{\alpha\in H^2(M)}S(\mathfrak{F}^\alpha)\hookrightarrow S(\mathfrak{F})\) is the standard inclusion. Here, \(\psi:\aut(K,\bar{\lambda})\rightarrow S(\mathfrak{F})\) is the natural projection.
\end{lemma}
\begin{proof}
  By Lemma \ref{sec:char-pairs-their-2}, a homomorphism \(\phi:\prod_{\alpha\in H^2(M)}S(\mathfrak{F}^\alpha)\hookrightarrow \aut(K,\bar{\lambda})\) with the properties described above is unique, if it exists.

  Therefore we only have to show that \(\phi\) exists.
  Let \(\sigma\in \prod_{\alpha\in H^2(M;\mathbb{Z})}S(\mathfrak{F}^\alpha)\subset S(\mathfrak{F})\).
  Then we have, for \(I\subset \mathfrak{F}\),
  \begin{equation*}
    \bigcap_{F_i\in I} F_i =\emptyset \;\Leftrightarrow \; \prod_{M_i\in I} PD(M_i)=\prod_{M_i\in I}PD(\sigma(M_i))=0 \; \Leftrightarrow \; \bigcap_{F_i\in I} \sigma(F_i)=\emptyset.
  \end{equation*}
Therefore \(\sigma\) is an automorphism of \(K\).

Now let \(F_1,\dots,F_n\in \mathfrak{F}\) such that \(\bigcap_{i=1}^n F_i\neq \emptyset\).
Then \(\bar{\lambda}(F_1),\dots,\bar{\lambda}(F_n)\) is a basis of \(IT\).
Moreover, since \(\sigma\) is an automorphism of \(K\), the same holds for  \(\bar{\lambda}(\sigma(F_1)),\dots,\bar{\lambda}(\sigma(F_n))\).

Therefore we may define an automorphism \(g_\sigma\) of \(T\) by \(Lg_\sigma(\bar{\lambda}(F_i))=\bar{\lambda}(\sigma(F_i))\) for \(i=1,\dots,n\).
Let \(v_1,\dots,v_n\) be the basis of \(IT^*\) dual to \(\bar{\lambda}(F_1),\dots,\bar{\lambda}(F_n)\) and \(v'_1,\dots,v'_n\) the basis of \(IT^*\) dual to \(\bar{\lambda}(\sigma(F_1)),\dots,\bar{\lambda}(\sigma(F_n))\).
Then we have, for \(i=1,\dots,n\),
\begin{align*}
  \sum_{M_j\in\mathfrak{F}-\{M_1,\dots,M_n\}} \langle v_i, \bar{\lambda}(M_j)\rangle PD(M_j)&= -PD(M_i)=-PD(\sigma(M_i))\\
&=\sum_{M_j\in\mathfrak{F}-\{\sigma(M_1),\dots,\sigma(M_n)\}} \langle v'_i, \bar{\lambda}(M_j)\rangle PD(M_j)\\
&=\sum_{M_j\in\mathfrak{F}-\{M_1,\dots, M_n\}} \langle v'_i, \bar{\lambda}(\sigma(M_j))\rangle PD(M_j).
\end{align*}
Since \(\{PD(M_i);\; M_i\in \mathfrak{F}-\{M_1,\dots,M_n\}\}\) is a basis of \(H^2(M)\), it follows that \(Lg_\sigma(\bar{\lambda}(M_i))=\bar{\lambda}(\sigma(M_i))\) for all \(M_i\in \mathfrak{F}\).

Therefore \((\sigma,g_\sigma)\in \aut(K,\bar{\lambda})\) and we define \(\phi(\sigma)=(\sigma,g_\sigma)\).
\end{proof}

\begin{lemma}
\label{sec:char-pairs-their-1}
  Let \(M\) be a quasitoric manifold with a strong omniorientation.
  Then for \(x\in M^T\) and \(\alpha\in H^2(M)\) we have
  \begin{equation*}
    \#\mathfrak{F}^\alpha-1\leq \#\{M_i\in \mathfrak{F}^\alpha;\; x\in M_i\} \leq \# \mathfrak{F}^\alpha.
  \end{equation*}
\end{lemma}
\begin{proof}
  It follows from the relations (\ref{eq:1}) that the Poincar\'e duals of the \(M_i\in \mathfrak{F}\), \(x\not\in M_i\), form a basis of \(H^2(M)\).
  Therefore there is at most one \(M_i\in \mathfrak{F}^\alpha\) which does not contain \(x\).
\end{proof}

\section{Constructing group actions}
\label{sec:construct}

In this section we construct an action of a compact connected Lie group on a quasitoric manifold which extends the torus action.

Before we do that, we explain how an action of a compact connected Lie group on a quasitoric manifold \(M\) induces a homomorphism of the Weyl group \(W(G)\rightarrow \aut(K,\bar{\lambda})\).
Here and in the following we choose a strong omniorientation of \(M\).

Assume that there is a compact connected Lie group \(G\) such that a finite quotient of \(G\) acts on \(M\) by an extension of the \(T\)-action.
Then \(M\) is called a quasitoric manifold with \(G\)-action.

The following structure results were shown in Section 2 of \cite{torus}. 
The group \(G\) has a covering group of the form \(\prod_{i=1}^kSU(l_i+1)\times T^{l_0}\).
Moreover, if \(g\in N_GT\) and \(M_i\in \mathfrak{F}\), then \(gM_i\) is a characteristic submanifold of \(M\).

Therefore we get an action of \(W(G)\) on \(\mathfrak{F}\).
Since \(G\) is connected, this action identifies \(W(G)\) with a subgroup of \(\prod_{\alpha\in H^2(M)}S(\mathfrak{F}^\alpha)\subset S(\mathfrak{F})\).
There are disjoint subsets \(\mathfrak{F}_1,\dots,\mathfrak{F}_k\) of \(\mathfrak{F}\) such that \(W(SU(l_i+1))\) is identified with \(S(\mathfrak{F}_i)\subset S(\mathfrak{F})\) for \(i=1,\dots,k\).
We call \(\mathfrak{F}_i\) the set of the characteristic submanifolds of \(M\) which are permuted by \(W(SU(l_i+1))\).

Moreover, \(W(G)\) acts on \(T\) by conjugation.
The actions of \(W(G)\) on \(T\) and \(\mathfrak{F}\) induce an action of \(W(G)\) on the characteristic pair \((K,\bar{\lambda})\).
The homomorphism \(W(G)\rightarrow \aut(K,\bar{\lambda})\) corresponding to this action is the restriction of the homomorphism \(\phi\) constructed in Lemma~\ref{sec:char-pairs-their} to \(W(G)\).

Now we state the main theorem of this section.

\begin{theorem}
\label{sec:constr-group-acti}
  Let \(M\) be a quasitoric manifold with a strong omniorientation. Then there is a smooth structure and a smooth action of a compact connected Lie group \(G\) on \(M\) which extends the torus action such that \(W(G)=\prod_{\alpha\in H^2(M)}S(\mathfrak{F}^\alpha)\).
\end{theorem}
\begin{proof}
  We prove this theorem by induction on the dimension of \(M\).
  If \(\dim M=0\) or \(\#\mathfrak{F}^\alpha\leq 1\) for all \(\alpha\in H^2(M)\), then there is nothing to prove.
  Therefore assume that there is an \(\alpha\in H^2(M)\) such that \(\#\mathfrak{F}^\alpha\geq 2\).
 
By Lemmas \ref{sec:char-pairs-their} and \ref{sec:char-pairs-their-1}, there are two cases:
\begin{enumerate}
\item \label{item:1}\(\bigcap_{M_i\in\mathfrak{F}^\alpha} M_i=\emptyset\) and, for all \(M_{i_0}\in \mathfrak{F}^\alpha\),  \(\bigcap_{M_i\in\mathfrak{F}^\alpha-\{M_{i_0}\}} M_i\neq\emptyset\),
\item \label{item:2} \(\bigcap_{M_i\in\mathfrak{F}^\alpha} M_i\neq\emptyset\).
\end{enumerate}

We consider at first the case (\ref{item:1}).
In this case we have with \(N=\bigcap_{M_i\in\mathfrak{F}^\alpha-\{M_{i_0}\}} M_i\) that
\begin{enumerate}
\item \label{item:3}\(M/T\) and \(\Delta^{\#\mathfrak{F}^\alpha-1}\times N/T\) are combinatorially equivalent and
\item\label{item:4}
  \begin{equation*}
    \Lambda(M)=
    \left(\begin{matrix}
      1& & &       -1& 0&\dots&0\\
       &\ddots & & \vdots& \vdots&\ddots&\vdots\\
       && 1&-1&0&\dots&0\\
      0 &\dots&0&\lambda_1&&&\\
      \vdots &\ddots&\vdots&\vdots&&\Lambda(N)&\\
      0 &\dots&0&\lambda_k&&&\\
    \end{matrix}\right),
  \end{equation*}
  where \(\Lambda(M)\) and \(\Lambda(N)\) are the characteristic matrices of \(M\) and \(N\), respectively. Here the first columns correspond to the facets in \(\mathfrak{F}^\alpha\).
They are of the form \(F\times N/T\), where \(F\) runs through the facets of \(\Delta^{\#\mathfrak{F}^\alpha-1}\).
\end{enumerate}

We show at first (\ref{item:3}). If this is shown then (\ref{item:4}) follows immediately from (\ref{eq:1}).
Let \(M_i\in \mathfrak{F}(M)-\mathfrak{F}^\alpha\), then by Lemmas \ref{sec:char-pairs-their} and \ref{sec:char-pairs-their-1} we must have \(N\cap M_i\neq \emptyset\).
Therefore we have a bijection \(\mathfrak{F}(M)-\mathfrak{F}^\alpha\rightarrow \mathfrak{F}(N)\).
Because, by the same lemmas, the intersection of \(M_{i_1},\dots,M_{i_k}\in \mathfrak{F}(M)-\mathfrak{F}^\alpha\) is empty if and only if \(N\cap \bigcap_{j=1}^kM_{i_j}=\emptyset\), (\ref{item:3}) follows.

Therefore, by Proposition 1.8 of \cite{0733.52006}, \(M\) is equivariantly homeomorphic to
\begin{equation*}
  S^{2\#\mathfrak{F}^\alpha -1}\times_{S^1}N.
\end{equation*}
Here the action of \(S^1\) on \(N\) is induced by the homomorphism \(\phi: S^1\rightarrow T'\) to the torus, which acts on \(N\), defined by \(\mu=(\lambda_1,\dots,\lambda_k)^t=\sum_{M_i\in\mathfrak{F}^\alpha}\bar{\lambda}(M_i)\).
Moreover, \(S^1\) acts on \(S^{2\#\mathfrak{F}^\alpha -1}\subset \C^{2\#\mathfrak{F}^\alpha}\) by multiplication.

By the induction hypothesis there is a compact connected Lie group \(G'\) which acts on \(N\) by an extension of the torus action, such that \(G'\) realizes the action of \(\prod_{\beta\in H^2(N)}S(\mathfrak{F}^\beta)\) on the simplicial complex dual to \(N/T\).

Then the action of \(SU(\#\mathfrak{F}^\alpha)\times G'\) on \(S^{2\#\mathfrak{F}^\alpha -1}\times N\) induces an action of \(G=SU(\#\mathfrak{F}^\alpha)\times Z_{G'}(\phi(S^1))\) on \(M\) such that the action of \(G/H\) extends the torus action.
Here \(Z_{G'}(\phi(S^1))\) is the centralizer of \(\phi(S^1)\) in \(G'\) and \(H\) is the ineffective kernel of the \(G\)-action on \(M\).
By Corollary 6.8 of \cite[p. 29]{MR0413144}, \(Z_{G'}(\phi(S^1))\) is connected.

From Lemma 2.5 of \cite{torus} we know that \(W(SU(\#\mathfrak{F}^\alpha))=S(\mathfrak{F}^\alpha)\).
Therefore we have to show that
\begin{equation*}
  W(Z_{G'}(\phi(S^1)))=\prod_{\beta\in H^2(M)-\{\alpha\}}S(\mathfrak{F}^\beta(M))\subset\prod_{\beta\in H^2(N)}S(\mathfrak{F}^\beta(N))=W(G').
\end{equation*}
It follows from the remarks at the beginning of this section that \(W(Z_{G'}(\phi(S^1)))\) is a subgroup of \(\prod_{\beta\in H^2(M)-\{\alpha\}}S(\mathfrak{F}^\beta(M))\).
Therefore let \(w\in\prod_{\beta\in H^2(M)-\{\alpha\}}S(\mathfrak{F}^\beta(M))\).
Since \(\prod_{\beta\in H^2(M)-\{\alpha\}}S(\mathfrak{F}^\beta(M))\) is generated by transpositions, we may assume that \(w\) is a transposition.

Let \(M_1,\dots,M_{m'}\in \mathfrak{F}(M)-\mathfrak{F}^\alpha\) such that \(\bigcap_{i=1}^{m'}M_i\cap N\) is a single point.
Then \(\bar{\lambda}(M_1),\dots,\bar{\lambda}(M_{m'})\) form a basis of \(IT'\).
Let \(v_1,\dots,v_{m'}\) be the dual basis of \(IT'^*\).
Let \(i\in\{1,\dots,m'\}\). At first assume that \(w(M_i)=M_{i'}\) with \(i'\in\{1,\dots,m'\}\). Then we have
\begin{align*}
  \langle v_i,\mu\rangle \alpha +\sum \langle v_i,\bar{\lambda}(M_j)\rangle PD(M_j)&=-PD(M_i)=-PD(M_{i'})\\
&=\langle v_{i'},\mu\rangle \alpha +\sum \langle v_{i'},\bar{\lambda}(M_j)\rangle PD(M_j).
\end{align*}
Here the sums are taken over those characteristic submanifolds of \(M\) which do not belong to \(\mathfrak{F}^\alpha\cup \{M_1,\dots,M_{m'}\}\).
Their Poincar\'e duals together with \(\alpha\) form a basis of \(H^2(M)\).
Therefore we have
\begin{equation*}
  \langle v_i,\mu\rangle=\langle v_{i'},\mu\rangle=\langle w^*v_{i},\mu\rangle =\langle v_{i},w_*\mu\rangle.
\end{equation*}

Now assume that \(w(M_i)\neq M_1,\dots,M_{m'}\). Then, because \(PD(M_i)=PD(w(M_i))\) we must have \(\langle v_i,\bar{\lambda}(M_j)\rangle =0\) for all \(M_j\in\mathfrak{F}(M)-\{ M_i,w(M_i)\}\). This implies \(\langle v_i,\mu\rangle=0\).
Moreover, we have
\begin{align*}
  \langle v_i,w_*\mu\rangle&=\sum_{j=1}^{m'}\langle v_i,w_*\bar{\lambda}(M_j)\rangle \langle v_j,\mu\rangle\\
&=\sum_{j\in\{1,\dots,m'\}-\{i\}}\langle v_i,\bar{\lambda}(M_j)\rangle \langle v_j,\mu\rangle=0.
\end{align*}

 Therefore we have \(\mu=w_*\mu\). This implies \(w\in W(Z_{G'}(\phi(S^1)))\).
 Hence, the claim follows in this case.

Now assume that \(\bigcap_{M_i\in \mathfrak{F}^\alpha}M_i\) is non-empty.
Then let \(\tilde{M}\) be the blow-up of \(M\) along  \(\bigcap_{M_i\in \mathfrak{F}^\alpha}M_i\) (see Section 4 of \cite{torus} for details).
If we write the characteristic matrix of \(M\) in the form
\begin{equation}
\label{eq:2}
  \left(
    \begin{matrix}
      1\\
      &\ddots&&\Lambda'\\
      &&1
    \end{matrix}\right)
\end{equation}
 such that the first \(\#\mathfrak{F}^\alpha\) columns correspond to the \(F_i\in \mathfrak{F}^\alpha\), then the characteristic matrix of \(\tilde{M}\) is given by
 \begin{equation*}
   \left(
     \begin{matrix}
       1&      & && &      & &&-1\\
        &\ddots& && &      & &&\vdots\\
        &      &1&& &      & &&-1\\
        &      & &1&&      & &\Lambda'&0\\
        &      & &&1&      & &&0\\
        &      & && &\ddots& &&\vdots\\
        &      & && &      &1&&0
     \end{matrix}
\right),
 \end{equation*}
where the proper transforms of the characteristic submanifolds of \(M\) are ordered as in (\ref{eq:2}), the last column corresponds to the exceptional submanifold and the first \(\#\mathfrak{F}^\alpha\) entries in this column are equal to \(-1\).

Hence two characteristic submanifolds of \(M\) have the same Poincar\'e duals if and only if their proper transforms have the same Poincar\'e duals.
Moreover, the Poincar\'e dual of the exceptional submanifold is distinct from the Poincar\'e duals of the other characteristic submanifolds of \(\tilde{M}\).

By the first case there is a \(G\)-action on \(\tilde{M}\) which extends the torus action.
Since the exceptional submanifold is fixed by \(\phi(S^1)\), we can \(G\)-equivariantly blow down \(\tilde{M}\) along the exceptional manifold to get a \(G\)-action on \(M\) (see \cite[Section 4]{torus} for details).
\end{proof}

\begin{cor}
  The group \(G\) constructed in Theorem \ref{sec:constr-group-acti} has a covering group of the form \(\prod_{\alpha\in H^2(M;\mathbb{Z})}SU(\#\mathfrak{F}^\alpha)\times T^{l_0}\).
\end{cor}
\begin{proof}
  This follows from the results of Section 2 of \cite{torus} and the description of the \(W(G)\)-action on \(\mathfrak{F}\) given in Theorem~\ref{sec:constr-group-acti}.
\end{proof}

In \cite{wiemeler_pre4} we proved that the equivariant smooth structures on a quasitoric manifold correspond one-to-one to the smooth structures on its orbit space.
We will show that the \(G\)-action constructed in Theorem~\ref{sec:constr-group-acti} is smooth with respect to the smooth structure on \(M\) which corresponds to the natural smooth structure on the simple polytope \(P\).
This will follow from the proof of Theorem~\ref{sec:constr-group-acti}, Corollary 5.3 of \cite{wiemeler_pre4} and the following lemma.

\begin{lemma}
\label{sec:smooth-struct-polyt-2}
   Let \(N\) be a quasitoric manifold and \(T^{l_0}\) the torus which acts on \(N\).
   Let, moreover, \(S^1\rightarrow T^{l_0}\) be a homomorphism and \(A\subset N\) be a characteristic submanifold, such that \(A\subset N^{S^1}\).
   We define \(\tilde{M}=S^{2k+1}\times_{S^1} N\) and \(M\) to be the blow-down of \(\tilde{M}\) along \(\mathbb{C}P^k\times A\).
   Then \(N/T^{l_0}\) is diffeomorphic to a simple polytope if and only if \(M/T\) is diffeomorphic to a simple polytope.
\end{lemma}
\begin{proof}
  It has been shown in the proof of Theorem 5.16 of \cite{torus} that \(N/T^{l_0}\) is homeomorphic to a simple polytope if and only if \(M/T\) is homeomorphic to a simple polytope.
  We will follow the proof of that theorem and show that all maps appearing there can be replaced by diffeomorphisms.
  Therefore the statement follows.

  Indeed, if \(M/T\) is diffeomorphic to a simple polytope, then \(N/T^{l_0}\) is diffeomorphic to a simple polytope because \(N/T^{l_0}\) is a face of \(M/T\).

  So we only have to prove the other implication.
  If \(N/T^{l_0}\) is diffeomorphic to a simple polytope, then all face-preserving homeomorphisms constructed in the proof of the cited theorem may be replaced by diffeomorphisms.
  So if we follow the proof of this theorem we end-up with a diffeomorphism
  \begin{equation*}
    g:F_1\times \Delta^{l_1}\rightarrow F_1\times \Delta^{l_1},
  \end{equation*}
which we want to extend to a diffeomorphism \(F_1\times \Delta^{l_1+1}\rightarrow F_1\times \Delta^{l_1+1}\), where \(F_1=A/T^{l_0}\).
Because every facet of \(F_1\times \Delta^{l_1}\) of the form \(F_1\times F\), where \(F\) is a facet of \(\Delta^{l_1}\), is mapped by \(g\) to a facet of the same form, it follows that the automorphism of the face poset of \(F_1\times \Delta^{l_1}\) induced by \(g\) is the restriction of an automorphism of the face poset of \(F_1\times \Delta^{l_1+1}\).
Therefore it follows from Theorem 5.1 of \cite{wiemeler_pre4} that \(g\) extends to a diffeomorphism of \(F_1\times \Delta^{l_1+1}\).
\end{proof}

\begin{cor}
  The action of the group \(G\) constructed in Theorem~\ref{sec:constr-group-acti} is smooth with respect to the smooth structure on \(M\) corresponding to the natural smooth structure on \(P\).
\end{cor}
\begin{proof}
  We use the same induction and notations as in the proof of Theorem~\ref{sec:constr-group-acti}.
  If the \(G'\)-action on \(N\) is smooth with respect to the smooth structure on \(N\) for which \(N/T\) is diffeomorphic to a simple polytope, then by Lemma~\ref{sec:smooth-struct-polyt-2} the \(G\)-action on \(M\) is smooth with respect to the smooth structure for which \(M/T\) is diffeomorphic to a simple polytope.
Since two simple polytopes are combinatorially equivalent if and only if they are diffeomorphic \cite[Corollary 5.3]{wiemeler_pre4}, it follows that the \(G\)-action on \(M\) is smooth with respect to the smooth structure for which \(M/T\) is diffeomorphic to \(P\).
\end{proof}

\section{Classification}
\label{sec:class}

Let \(G\) be a compact Lie group with maximal torus \(T\).
A quasitoric manifold with \(G\)-action is a smooth \(G\)-manifold \(M\) such that \(M\) together with the action of the maximal torus of \(G/H\) is a quasitoric manifold.
Here \(H\) is a finite subgroup of \(G\) which acts trivially on \(M\).
Examples of such manifolds are \(\C P^{l_1}\) or bundles with fiber a \(2l_0\)-dimensional quasitoric manifold and structure group \(S^1\) over \(\C P^{l_1}\).
In these cases we have \(G=SU(l_1+1)\) or \(G=SU(l_1+1) \times T^{l_0}\), respectively.

In \cite{torus} we classified quasitoric manifolds with \(G\)-action.
As a first step towards this classification we showed that \(G\) has a covering group of the form \(\prod_{i=1}^kSU(l_i+1)\times T^{l_0}\).

The classification was given in terms of admissible triples \((\psi,N,(A_1,\dots,A_k))\), where
\begin{itemize}
\item \(\psi\) is an homomorphism \(\prod_{i=1}^kS(U(l_i)\times U(1))\rightarrow T^{l_0}\).
\item \(N\) is a \(2l_0\)-dimensional quasitoric manifold.
\item The \(A_i\) are pairwise distinct characteristic submanifolds of \(N\) or empty.
  If \(A_i\) is non-empty then \(\image \psi|_{S(U(l_i)\times U(1))}\) acts trivially on \(A_i\) and
  \begin{equation*}
    \ker \psi|_{S(U(l_i)\times U(1))}=SU(l_i).
  \end{equation*}
\end{itemize}

Two such triples \((\psi,N,(A_1,\dots,A_k))\) and \((\psi',N',(A_1',\dots,A_k'))\) are called equivalent or diffeomorphic if
\begin{itemize}
\item \(\psi|_{S(U(l_i)\times U(1))}=\psi'|_{S(U(l_i)\times U(1))}\) if \(l_i>1\).
\item \(\psi|_{S(U(l_i)\times U(1))}=\psi'^{\pm 1}|_{S(U(l_i)\times U(1))}\) if \(l_i=1\).
\item There is an \(T^{l_0}\)-equivariant diffeomorphism \(f:N\rightarrow N'\) such that \(f(A_i)=A'_i\) for all \(i\).
\end{itemize}
The main theorem of \cite{torus} may be formulated in the following way:
\begin{theorem}[{\cite[Theorem 8.6]{torus}}]
\label{sec:classification-3}
  Let \(G=\prod_{i=1}^k SU(l_i+1)\times T^{l_0}\). Then the \(G\)-equivariant diffeomorphism classes of quasitoric manifolds with \(G\)-action are in one-to-one correspondence with the diffeomorphism classes of admissible triples.
\end{theorem}

If in the situation of the above theorem \(G\) has only one simple factor, then the one-to-one correspondence is given by the following construction (for details see Section 5 of \cite{torus}).
Let \(M\) be a quasitoric manifold with \(G\)-action and \((\psi,N,A_1)\) the admissible triple corresponding to \(M\).
Then \(M\) can be reconstructed from \((\psi,N,A_1)\) as follows.

Let \(\tilde{M}=SU(l_1+1)\times_{S(U(l_1)\times U(1))}N\), where the \(S(U(l_1)\times U(1))\)-action on \(N\) is induced by \(\psi^{-1}\).
Then \(M\) is the blow-down of \(\tilde{M}\) along \(SU(l_1+1)/S(U(l_1)\times U(1))\times A_1\).
It follows from these constructions that there is a natural identification \(M/G=N/T^{l_0}\) (see the proof of Corollary 8.8 of \cite{torus} for details).
Under this identification \(M^{SU(l_1+1)}/G\) is identified with \(A_1/T^{l_0}\).

In the other direction \((\psi,N,A_1)\) is determined by \(M\) as follows.
There is a characteristic submanifold \(M_{i_0}\in \mathfrak{F}_1\) which is fixed by the action of \(W(S(U(l_1)\times U(1)))\) on \(\mathfrak{F}\).
\(N\) is the intersection of the other characteristic submanifolds belonging to \(\mathfrak{F}_1\).
Moreover, \(A_1=M^{SU(l_1+1)}=\bigcap_{M_i\in\mathfrak{F}_1} M_i\).

The homomorphism \(\psi\) can be defined as follows.
Let \(x\in N^T\). Then \(x\) is also fixed by \(S(U(l_1)\times U(1))\) and \(\psi:S(U(l_1)\times U(1))\rightarrow T^{l_0}\) is the unique homomorphism such that all \((g,\psi(g))\), \(g\in S(U(l_1)\times U(1))\), act trivially on \(T_xN\) and \(N\).

If there are more than one simple factors, say \(G=\prod_{i=1}^k SU(l_i+1)\times T^{l_0}\), then
the submanifold \(N\) from above is invariant under the action of \(G'=\prod_{i=2}^k SU(l_i+1)\times T^{l_0}\).
Hence it is a quasitoric manifold with \(G'\)-action.
So one can iterate the above constructions to get a classification of quasitoric manifolds with \(G\)-action (see Section 8 of \cite{torus} for details).

We call two admissible triples \((\psi,N,(A_1,\dots,A_k))\) and \((\psi',N',(A_1',\dots,A_k'))\) homeomorphic if
\begin{itemize}
\item \(\psi|_{S(U(l_i)\times U(1))}=\psi'|_{S(U(l_i)\times U(1))}\) if \(l_i>1\).
\item \(\psi|_{S(U(l_i)\times U(1))}=\psi'^{\pm 1}|_{S(U(l_i)\times U(1))}\) if \(l_i=1\).
\item There is a \(T^{l_0}\)-equivariant homeomorphism \(f:N\rightarrow N'\) such that \(f(A_i)=A'_i\) for all \(i\).
\end{itemize}

With this notation we have the following classification of quasitoric manifolds with \(G\)-action up to \(G\)-equivariant homeomorphism.

\begin{theorem}
\label{sec:classification}
  Let \(G=\prod_{i=1}^k SU(l_i+1)\times T^{l_0}\). Then the \(G\)-equivariant homeomorphism classes of quasitoric manifolds with \(G\)-action are in one-to-one correspondence with the homeomorphism classes of admissible triples.
\end{theorem}
\begin{proof}
  At first note that the homeomorphism type of the admissible triple corresponding to a quasitoric manifold with \(G\)-action depends only on the \(G\)-equivariant homeomorphism type of \(M\) because \(N\) and the \(A_i\) may be identified with intersections of characteristic submanifolds of \(M\) and \(\psi\) depends only on the isotropy groups of points in these intersections.

Because as noted above \(M/G=N/T^{l_0}\), it follows from Lemma~\ref{sec:classification-2} below that the homeomorphism type of the admissible triple determines the \(G\)-equivariant homeomorphism type of \(M\). 
\end{proof}

By the results of Section 5 of \cite{wiemeler_pre4}, two quasitoric manifolds of dimension at most \(6\) are equivariantly diffeomorphic if and only if they are equivariantly homeomorphic.
Hence, Theorems \ref{sec:classification-3} and  \ref{sec:classification} show that quasitoric manifolds with \(G\)-action, where the center of \(G\) has dimension at most \(3\), are equivairantly homeomorphic if and only if they are equivariantly diffeomorphic.
This might fail if the dimension of the center is higher.

\begin{lemma}
\label{sec:classification-2}
  Let \(M\) be a quasitoric manifold with \(G\)-action, \(G=\prod_{i=1}^k SU(l_i+1)\times T^{l_0}\).
Then \(M\) is equivariantly homeomorphic to
\( M/G\times G/\sim \), where \((x,g)\sim(x',g')\) if and only if \(x=x'\) and \(g^{-1}g'\in H_x\).
Here the groups \(H_x\) depend only on \(x\) and the admissible triple corresponding to \(M\).
\end{lemma}
\begin{proof}
  We prove this lemma by induction on the number of simple factors of \(G\).
If \(G\) is a torus then this lemma is due to Davis and Januszkiewicz \cite[Proposition 1.8]{0733.52006}.

Therefore we may assume that there is at least one simple factor.
There are two cases:
\begin{enumerate}
\item \(M^{SU(l_1+1)}=\emptyset\)
\item \(M^{SU(l_1+1)}\neq \emptyset\).
\end{enumerate}

In the first case we have by Corollary 5.6 of \cite{torus}
\begin{equation*}
  M=SU(l_1+1)\times_{S(U(l_1)\times U(1))}M', 
\end{equation*}
where \(M'\) is a quasitoric manifold with \(G'\)-action, \(G'=\prod_{i=2}^kSU(l_i+1)\times T^{l_0}\).
The action of \(S(U(l_1)\times U(1))\) on \(M'\) is induced by the homomorphism
\begin{equation*}
   \phi=(\psi|_{S(U(l_i)\times U(1))})^{-1}:S(U(l_1)\times U(1))\rightarrow T^{l_0},
\end{equation*}
 where \(\psi\) is the homomorphism from the admissible triple of \(M\).
Since \(M/G=M'/G'\), we have
by the induction hypothesis
\begin{equation*}
  M= SU(l_1+1)\times_{S(U(l_1)\times U(1))} \left(M/G\times G'/\sim'\right).
\end{equation*}
For a subgroup \(H'_x\) of \(G'\) we have
\begin{equation*}
  SU(l_1+1)\times_{S(U(l_1)\times U(1))}(G'/H'_x)= (SU(l_1+1)\times G')/(\phi \id_{G'})^{-1}(H'_x).
\end{equation*}
Therefore the statement follows in this case with \(H_x=(\phi \id_{G'})^{-1}(H'_x)\).

If \(M^{SU(l_1+1)}\) is non-empty.
Then \(M\) is the blow-down of some quasitoric manifold \(\tilde{M}\) with \(G\)-action but without \(SU(l_1+1)\)-fixed points. Let \(F:\tilde{M}\rightarrow M\) be the projection.
Then we have \(M=\tilde{M}/\sim''\), where \(y\sim''y'\) if and only if there is a \(g\in SU(l_1+1)\) such that \(gy=y'\) in the case \(y,y'\in F^{-1}(M^{SU(l_1+1)})\) or \(y=y'\) otherwise.

Under the identification of \(M/G\) with \(N/T^{l_0}\) given above, \(M^{SU(l_1+1)}/G\) is identified with \(A_1/T^{l_0}\).
Hence, for \(x \in M^{SU(l_1+1)}/G\), we have \(\image \phi=\image\psi\subset H'_x\).

Therefore we have \(S(U(l_1)\times U(1))\times H'_x=(\phi \id_{G'})^{-1}(H'_x)\) if \(x\in M^{SU(l_1+1)}/G\).
Hence, the statement follows with
\begin{equation*}
  H_x=
  \begin{cases}
    SU(l_1+1)\times H'_x&\text{if } x\in M^{SU(l_1+1)}/G\\
    (\phi \id_{G'})^{-1}(H'_x)&\text{if } x\not\in M^{SU(l_1+1)}/G.
  \end{cases}
\end{equation*}
\end{proof}

\begin{theorem}
\label{sec:classification-1}
  Let \(G=\prod_{i=1}^kSU(l_i+1)\times T^{l_0}\) and \(M\), \(M'\) be two quasitoric manifolds with \(G\)-action. Let \(T\) be a maximal torus of \(G\).
  Then \(M\) and \(M'\) are \(G\)-equivariantly homeomorphic (diffeomorphic), if and only if they are \(T\)-equivariantly homeomorphic (diffeomorphic).
\end{theorem}
\begin{proof}
  Without loss of generality we may assume that \(T\) is the standard maximal torus of \(G\).
  Let \((\psi,N,(A_1,\dots,A_k))\) be the admissible triple corresponding to \(M\).
  We show that \((\psi,N,(A_1,\dots,A_k))\) is determined up to homeomorphism (diffeomorphism) by the homeomorphism (diffeomorphism) type of the \(T\)-action on \(M\).

Let \(\pi_i:G\rightarrow SU(l_i+1)\) the projection and \(\mathfrak{F_i}\) the set of characteristic submanifolds which are permuted by \(W(SU(l_i+1))\).
By Lemma 2.10 of \cite{torus}, we have, for \(w\in W(SU(l_i+1))\) and \(M_j\in \mathfrak{F}\), \(w\lambda(M_j)w^{-1}=\lambda(wM_j)\).
Therefore \(\pi_i\circ \lambda(M_j)\) is trivial if \(M_j\not\in \mathfrak{F}_i\).
Hence, there is an \(M_{j_0}\in \mathfrak{F}_i\) such that \(\pi_i\circ \lambda(M_{j_0})\) is non-trivial.
By Lemma 2.7 of \cite{torus}, \(W(SU(l_i+1))\) acts transitively on \(\mathfrak{F}_i\).
Hence, \(\pi_i\circ\lambda(M_j)\) is non-trivial if and only if  \(M_j\in \mathfrak{F}_i\).

If \(l_i>1\), there is exactly one \(M_{j_i}\in \mathfrak{F}_i\) such that \(\lambda(M_{j_i})\) is fixed by \(W(S(U(l_i)\times U(1)))\).
Since each \(T\)-fixed point is contained in at least \(l_i\) characteristic submanifolds belonging to \(\mathfrak{F}_i\), we have, for each \(M_{j_0}\in \mathfrak{F}_i\), \(\dim \langle \pi_i\circ \lambda(M_j);\;M_j\in \mathfrak{F}_i -\{M_{j_0}\}\rangle=l_i\).
Because the center of \(S(U(l_1)\times U(1))\) is one-dimensional, \(M_{j_i}\) is the only characteristic submanifold such that \(\pi_i\circ \lambda(M_{j_i})\) is contained in the center of \(S(U(l_i)\times U(1))\).

If \(l_i=1\), then \(\mathfrak{F}_i\) has exactly two elements and any choice of a \(M_{j_i}\in \mathfrak{F}_i\) leads to the same equivalence class of admissible triples (see \cite[Section 5]{torus} for details).

Then we have
\begin{align*}
  N&=\bigcap_{i=1}^k\bigcap_{M_j\in\mathfrak{F}_i-\{M_{j_i}\}}M_j&
  \text{and}&&A_i&=N\cap M_{j_i}.
\end{align*}

By its construction in the proof of Lemma 5.3 of \cite{torus}, the homomorphism \(\psi\) depends only on the \(\prod_{i=1}^kS(U(l_i)\times U(1))\times T^{l_0}\)-representation \(T_xM\) with \(x\in N^{T^{l_0}}\).
Since \(T\) is a maximal torus of \(\prod_{i=1}^kS(U(l_i)\times U(1))\times T^{l_0}\), this representation depends only on the \(T\)-equivariant homeomorphism type of \(M\).

Now the statement follows from Theorem~\ref{sec:classification}.
\end{proof}

\section{Uniqueness}
\label{sec:unique}

In this section we prove that the group constructed in Section \ref{sec:construct} is a maximal compact connected Lie subgroup of the homeomorphism group of \(M\) which contains the torus and that it is unique up to conjugation.

\begin{lemma}
  \label{sec:uniqueness-1}
  Let \(M\) be a quasitoric manifold with \(G\)-action, \(G= \prod^{k}_{i=1} SU(l_i+1)\times T^{l_0}\) and \(T\) a maximal torus of \(G\), such that \(T^{l_0}\) acts effectively on \(M\).
  Denote by \(\mathfrak{F}_i\), \(i=1,\dots,k\), the set of characteristic submanifolds of \(M\) which are permuted by \(W(SU(l_i+1))\).
  Moreover, let  \(\mathfrak{F}_0=\mathfrak{F}-\bigcup_{i=1}^k\mathfrak{F}_i\).
  Then we have:
  \begin{enumerate}
  \item \label{item:5} The subgroup of \(G\) which acts trivially on \(M\) is given by 
    \begin{equation*}
      H=\{(g,\psi(g))\in G;\; g\in Z(\prod_{i=1}^kSU(l_i+1))\}.
    \end{equation*}
  \item \label{item:6} Let \(M_1\) be a characteristic submanifold of \(M\) which belongs to \(\mathfrak{F}_i\) and \(x\in M_1\) a generic point. Denote the identity component of \(T_x\) by \(T^0_x\).
    Then we have:
    \begin{equation*}
      H\cap T_x^0=\{(g,\psi(g))\in G;\; g\in Z(SU(l_i+1))\}
    \end{equation*}
    if \(i>0\). If \(i=0\) then \(H\cap T_x^0=1\).
  \end{enumerate}
\end{lemma}
\begin{proof}
  At first we prove (\ref{item:5}). We prove this statement by induction on \(k\).
  If \(k=0\), then there is nothing to prove.
  Therefore assume that \(k>0\) and that the statement is proved for all quasitoric manifolds with \(G'=\prod_{i=2}^kSU(l_i+1)\times T^{l_0}\)-action.
  With the notation from the proof of Lemma~\ref{sec:classification-2} the subgroup of \(G\) which acts trivially on \(M\) is given by
  \begin{align*}
    H&=\bigcap_{x\in M/G, g\in G}gH_xg^{-1}=\bigcap_{x\in M/G, g\in G} g(\phi\id_{G'})^{-1}(H_x')g^{-1}\\
    &=\bigcap_{g\in SU(l_1+1)}g(\phi\id_{G'})^{-1}(\bigcap_{x\in M/G, g'\in G'}g'H_x'g'^{-1})g^{-1}\\
    &=\bigcap_{g\in SU(l_1+1)}g\langle\{(h,\psi(h));\; h\in S(U(l_1)\times U(1))\}, \bigcap_{x\in M/G, g'\in G'}g'H_x'g'^{-1}\rangle g^{-1}\\
    &=\langle \{(h,\psi(h);\; h\in Z(SU(l_1+1))\}, H'\rangle\\
    &=\{(g,\psi(g))\in G;\; g\in Z(\prod_{i=1}^kSU(l_i+1))\}.
  \end{align*}
Here \(H'\) denotes the subgroup of \(G'\) which acts trivially on \(N\).

Now we prove the second statement.
At first assume \(i>0\).
After blowing up \(M\) along the fixed points of \(SU(l_i+1)\), we may assume that
\begin{equation*}
  M=SU(l_i+1)\times_{S(U(l_i)\times U(1))}N.
\end{equation*}
Then there is an \(SU(l_i+1)\)-equivariant projection \(p:M\rightarrow \C P^{l_i}\).
The characteristic submanifold \(M_1\) of \(M\) is given by a preimage of a characteristic submanifold \(\C P^{l_i}_1\) of \(\C P^{l_i}\).
Now we have
\begin{equation*}
  T_x^0=\{(t,\psi(t))\in T;\; t\in T^{l_i}_{p(x)}\},
\end{equation*}
where \(T^{l_i}\) denotes \(T\cap SU(l_i+1)\).
Since \(T^{l_i}_{p(x)}\) contains the center of \(SU(l_i+1)\) the statement follows in this case.

If \(i=0\), then it follows from Lemma 2.10 of \cite{torus} that \(T_x^0\) is fixed pointwise by the action of \(W(G)\) on \(T\).
Hence, it is contained in \(T^{l_0}\).
Therefore the statement follows in this case.
\end{proof}

\begin{lemma}
\label{sec:uniqueness}
  Let \(M\) be a quasitoric manifold and \(T\) the torus which acts on M by \(\phi:T\rightarrow \Homeo(M)\).
  Let \(G_j\), \(j=1,2\), be compact connected Lie groups and \(\iota_j:T\rightarrow G\), \(j=1,2\), embeddings of \(T\) as maximal tori of \(G_j\).
  Assume that there are effective actions \(\phi_j:G_j\rightarrow \Homeo(M)\), \(j=1,2\), such that \(\phi=\phi_j\circ \iota_j\) for \(j=1,2\).
  Moreover, assume that the natural actions of \(W(G_j)\), \(j=1,2\), on \(\mathfrak{F}\) induce identifications of \(W(G_1)\) and \(W(G_2)\) with a given subgroup \(H\) of \(S(\mathfrak{F})\).
Then \(\phi_1(G_1)\) and \(\phi_2(G_2)\) are conjugated in \(\Homeo(M)\).
\end{lemma}
\begin{proof}
  Since the Weyl groups of \(G_1\) and \(G_2\) are isomorphic.
  There is a group of the form \(\tilde{G}=\prod^{k}_{i=1}SU(l_i+1)\times T^{l_0}\) and coverings \(\varphi_j:\tilde{G}\rightarrow G_j\), \(j=1,2\).

  Because all maximal tori in \(\tilde{G}\) are conjugated, we may assume that there is a maximal torus \(\tilde{T}\) of \(\tilde{G}\) such that \(\tilde{T}=\varphi_j^{-1}(\iota_j(T))\), \(j=1,2\).
  Let \(\psi\) be the automorphism of \(L\tilde{T}\) given by \((L\varphi_1)^{-1}\circ L\iota_1\circ (L\iota_2)^{-1} \circ L\varphi_2\).
  
  Equip \(M\) with a strong omniorientation.
This omniorientation is preserved by the actions of \(G_1\) and \(G_2\).
Denote by \(\bar{\lambda}\) the characteristic function for the \(T\)-action on \(M\).
  Moreover, denote by \(\bar{\lambda}_j(M_i)\), \(j=1,2\), \(M_i\in \mathfrak{F}\), a primitive vector in \(I\tilde{T}\) which generates the identity component of the isotropy group of a generic point in \(M_i\) with respect to the \(\tilde{T}\)-action \(\phi_j\circ \varphi_j\).
We choose this primitive vector in such a way that it is compatible with the omniorientation chosen above. 
  
  Then, by Lemma \ref{sec:uniqueness-1}, we have, for all \(M_k\in \mathfrak{F}_i\), \(i>0\), and \(j=1,2\), 
  \begin{equation*}
    L\iota_j^{-1}\circ L\varphi_j(\bar{\lambda}_j(M_k))=(l_i+1)\bar{\lambda}(M_k).
  \end{equation*}
For \(M_k\in \mathfrak{F}_0\) we have
\begin{equation*}
   L\iota_j^{-1}\circ L\varphi_j(\bar{\lambda}_j(M_k))=\bar{\lambda}(M_k).
\end{equation*}

This implies that \(\psi(\bar{\lambda}_2(M_k))=\bar{\lambda}_1(M_k)\).
By Lemma 2.10 of \cite{torus}, we have, for \(w\in W(\tilde{G})\) and \(i=1,2\), \(\bar{\lambda}_i(wM_k)=w\bar{\lambda}_i(M_k)w^{-1}\).
Hence, it follows that
  \begin{equation*}
     \psi(w\bar{\lambda}_2(M_k)w^{-1})=\psi(\bar{\lambda}_2(wM_k))=\bar{\lambda}_1(wM_k)=w\bar{\lambda}_1(M_k)w^{-1}=w\psi(\bar{\lambda}_2(M_k))w^{-1}.
  \end{equation*}
It follows that \(\psi\) is an automorphism of the \(W(\tilde{G})\)-representation \(L\tilde{T}\).
Because each irreducible non-trivial summand of \(L\tilde{T}\) appears only once in a decomposition of \(L\tilde{T}\) in irreducible representations, it follows from Schur's Lemma that the restriction of \(\psi\) to the Lie-algebra of the maximal torus \(\tilde{T}_i\) of a simple factor \(SU(l_i+1)\) of \(\tilde{G}\) is multiplication with a constant \(a_i\in \R\).

Therefore we have 
\begin{equation*}
\iota_1^{-1}\circ \varphi_1(\tilde{T}_i)=\iota_2^{-1}\circ \varphi_2(\tilde{T}_i).
\end{equation*}
Denote this subtorus of \(T\) by \(T_i\).
By Lemma~\ref{sec:uniqueness-1}, we have that
\begin{equation*}
  IT_i/L\iota_1^{-1}\circ L\varphi_1(I\tilde{T}_i)\cong \ker \iota_1^{-1}\circ \varphi_1\cap \tilde{T}_i \cong \ker \iota_2^{-1}\circ \varphi_2\cap \tilde{T}_i \cong IT_i/L\iota_2^{-1}\circ L\varphi_2(I\tilde{T}_i).
\end{equation*}
Note that \(I_{ij}=\langle \bar{\lambda}_j(M_k);\; M_k\in \mathfrak{F}\rangle\cap L\tilde{T}_i\) is a lattice of maximal rank in \(I\tilde{T}_i\).
Then we have
\begin{align*}
  |I\tilde{T}_i/I_{i1}|&=\frac{|IT_i/L\iota_1^{-1}\circ L\varphi_1(I_{i1})|}{|IT_i/L\iota_1^{-1}\circ L\varphi_1(I\tilde{T}_i)|} 
  =\frac{|IT_i/L\iota_2^{-1}\circ L\varphi_2(I_{i2})|}{|IT_i/L\iota_2^{-1}\circ L\varphi_2(I\tilde{T}_i)|}\\
  &=|I\tilde{T}_i/I_{i2}|=|\psi(I\tilde{T}_i)/\psi(I_{i2})|=\frac{1}{|a_i|^{l_i}}|I\tilde{T}_i/I_{i1}|,
\end{align*}
because \(\psi(I_{i2})=I_{i1}\).
Therefore we must have \(a_i=\pm 1\).
Therefore there is an automorphism \(\Psi\) of \(\tilde{G}\) with \(L\Psi = \psi\).
Now the statement follows from Theorem \ref{sec:classification-1} applied to the \(\tilde{G}\)-actions \(\phi_1\circ\varphi_1\circ\Psi\) and \(\phi_2\circ \varphi_2\).
\end{proof}

\begin{remark}
  If, in the situation of Lemma~\ref{sec:uniqueness}, both \(G_j\)-actions are smooth with respect to the same smooth structure, then it follows from Theorem \ref{sec:classification-1} that \(\phi_1(G_1)\) and \(\phi_2(G_2)\) are conjugated in \(\Diff(M)\).
\end{remark}

\begin{theorem}
  Let \(M\) be a quasitoric manifold and \(G\subset \Homeo(M)\), the group constructed in Section~\ref{sec:construct}. If \(G'\) is another compact connected Lie group which acts by an extension of the torus action on \(M\). Then \(G'\) is conjugated in \(\Homeo(M)\) to a subgroup of \(G\).
If the \(G\) and \(G'\)-actions are smooth for the same smooth structure on \(M\), then \(G'\) is conjugated in \(\Diff(M)\) to a subgroup of \(G\).
\end{theorem}
\begin{proof}
Let \(\mathfrak{F}=\mathfrak{F}_1\amalg\dots\amalg\mathfrak{F}_k\) be a partition in \(W(G')\)-orbits. Then we have \(W(G')=\prod_{i=1}^kS(\mathfrak{F}_k)\).
Moreover, since the \(G'\)-action on \(H^*(M)\) is trivial, it follows that the sets \(\mathfrak{F}^\alpha\), \(\alpha\in H^2(M)\), are \(W(G')\)-invariant.

This gives as an homomorphism \(W(G')\rightarrow W(G)=\prod_{\alpha\in H^2(M)}S(\mathfrak{F}^\alpha)\).
There is a subgroup of maximal rank of \(G\) whose Weyl group is given by the image of this homomorphism.
Therefore the statement follows from Lemma \ref{sec:uniqueness}.
\end{proof}

Now we have proven all parts of Theorem \ref{sec:introduction} besides the statement about the symplectic toric manifolds.
To prove this part we first recall the construction of a maximal compact Lie subgroup of the symplectomorphism group of a symplectic toric manifold due to Masuda \cite{1210.53081}.
An alternative construction of this group was given by McDuff and Tolman \cite{1202.52010}.

Masuda showed that there is a root system \(R(M)\) such that the root system of every compact connected Lie subgroup of the symplectomorphism group which contains the torus is a subroot system of \(R(M)\).
Moreover, he constructed a compact Lie subgroup \(G'\) of the symplectomorphism group which contains the torus and has a root system isomorphic to \(R(M)\).

The proof of the first part of Masuda's results is also valid for any compact connected Lie subgroup of the homeomorphism group of \(M\) which contains the torus and  preserves the omniorientation induced by the symplectic form on \(M\).
Therefore \(G\) and \(G'\) are conjugated, if the \(G\)-action preserves this omniorientation.
Hence, it is sufficient to prove that this omniorientation is strong.

Let \(M_1\) and \(M_2\) be two characteristic submanifolds such that
\begin{equation*}
   PD(M_1)=\pm PD(M_2).
\end{equation*}

Let \(e_1,\dots,e_n\) be the standard basis of \(\R^n\).
At first assume that \(M_1\cap M_2=\emptyset\).
In this case we may assume that \(\{0\}=F_1\cap \bigcap_{i=3}^{n+1} F_i\subset \R^n\) and \(\bar{\lambda}(F_1)=e_1\), \(\bar{\lambda}(F_i)=e_{i-1}\) for \(i=3,\dots,n+1\).
It follows from \(PD(M_1)=\pm PD(M_2)\) that \(\bar{\lambda}(F_2)=\pm e_1 + \sum_{i=2}^{n}\mu_{i2}e_i\) and \(\bar{\lambda}(F_j)=\sum_{i=2}^n\mu_{ij}e_i\) for \(j>n+1\) with \(\mu_{ij}\in \mathbb{Z}\). 
We should note that if \(M\) is a symplectic toric manifold then \(\bar{\lambda}(F_i)\) is the outward normal vector of the facet \(F_i\) of \(P\).
Therefore \(P\cap \langle e_1 \rangle\) is an interval with boundary \(\langle e_1\rangle \cap (F_1\cup F_2)\).
Hence we must have  \(\bar{\lambda}(F_2)=-e_1+ \sum_{i=2}^{n}\mu_{i2}e_i\).
This implies \(PD(M_1)=PD(M_2)\).

Now consider the case \(M_1\cap M_2\neq\emptyset\).

Without loss of generality we may assume that \(\{0\}=\bigcap_{i=1}^n F_i\subset \R^n\) and \(\bar{\lambda}(F_i)=e_i\) for \(i=1,\dots n\).

Assume that \(PD(M_1)=-PD(M_2)\). 
Then for all \(F_j\in \mathfrak{F}\), \(j>n\), there are \(\mu_{0j},\mu_{3j},\dots,\mu_{nj}\in\mathbb{Z}\) such that
\begin{equation*}
  \bar{\lambda}(F_j)=\mu_{0j}(e_1-e_2)+\sum_{i=3}^n \mu_{ij}e_i.
\end{equation*}
Because the \(\bar{\lambda}(F_j)\) are the outward normal vectors of the facets of \(P\) it follows that \(P\cap \langle e_1,e_2\rangle\) is non-compact.
But this is impossible because \(P\) is a convex polytope.
Therefore we must have \(PD(M_1)=PD(M_2)\).

\bibliography{diff_qt}{}

\providecommand{\bysame}{\leavevmode\hbox to3em{\hrulefill}\thinspace}
\providecommand{\MR}{\relax\ifhmode\unskip\space\fi MR }
\providecommand{\MRhref}[2]{%
  \href{http://www.ams.org/mathscinet-getitem?mr=#1}{#2}
}
\providecommand{\href}[2]{#2}
\begin{thebibliography}{1}

\bibitem{MR0413144}
Glen~E. Bredon, \emph{Introduction to compact transformation groups}, Academic
  Press, New York, 1972, Pure and Applied Mathematics, Vol. 46. \MR{0413144 (54
  \#1265)}

\bibitem{0733.52006}
M.~W. Davis and T.~Januszkiewicz, \emph{{Convex polytopes, Coxeter orbifolds
  and torus actions.}}, Duke Math. J. \textbf{62} (1991), no.~2, 417--451
  (English).

\bibitem{1210.53081}
M.~Masuda, \emph{{Symmetry of a symplectic toric manifold.}}, J. Symplectic
  Geom. \textbf{8} (2010), no.~4, 359--380 (English).

\bibitem{1202.52010}
D.~McDuff and S.~Tolman, \emph{{Polytopes with mass linear functions. I.}},
  Int. Math. Res. Not. \textbf{2010} (2010), no.~8, 1506--1574 (English).

\bibitem{torus}
M.~Wiemeler, \emph{Torus manifolds with non-abelian symmetries}, Trans. Am.
  Math. Soc. \textbf{364} (2012), no.~3, 1427--1487.

\bibitem{wiemeler_pre4}
\bysame, \emph{Exotic torus manifolds and equivariant smooth structures on
  quasitoric manifolds}, Math. Z. \textbf{273} (2013), 1063--1084.

\end{thebibliography}
\bibliographystyle{amsplain}
\end{document}